\documentclass[10pt]{article}
\usepackage{latexsym,amsmath,bez123,calc,curves,ebezier,epic,eepic,graphicx,multiply,rotating}
\everymath{\displaystyle}
\usepackage{graphicx}
\usepackage{pst-all}
\usepackage{color}
\usepackage{amssymb}
\usepackage{amsmath}
\usepackage{pstricks}
\usepackage{epsfig}
\usepackage{pst-grad} 
\usepackage{pst-plot} 
\newtheorem{prethm}{{\bf Theorem}}

\newenvironment{thm}{\begin{prethm}{\hspace{-0.5
               em}{\bf .}}}{\end{prethm}}
\newtheorem{prelemma}{{\bf Lemma}}

\newenvironment{lemma}{\begin{prelemma}{\hspace{-0.5
               em}{\bf .}}}{\end{prelemma}}
\newtheorem{preex}{{\bf Example}}

\newtheorem{preprop}{{\bf Proposition}}

\newenvironment{prop}{\begin{preprop}{\hspace{-0.5em}{\bf .}}}{\end{preprop}}
\newtheorem{precor}{{\bf Corollary}}

\newenvironment{cor}{\begin{precor}{\hspace{-0.5
               em}{\bf .}}}{\end{precor}}
\newtheorem{preremark}{{\bf Remark}}

\newtheorem{preprob}{{\bf Problem}}

\newtheorem{predefin}{{\bf Definition}}

\newenvironment{defin}{\begin{predefin}{\hspace{-0.5
               em}{\bf .}}}{\end{predefin}}
\newtheorem{preconj}{{\bf Conjecture}}

\newtheorem{preprobb}{{\bf Problem}}

\newtheorem{prelem}{{\bf Theorem}}

\newenvironment{proof}{{\bf Proof.}\rm }{\hfill{$\Box$}}

\newtheorem{presolution}{{\bf Solution.}}

\def\newpic#1{}

\title{\Large\bf\noindent Weak dynamic monopolies in social graphs}

\author{\large\bf Mitra Nemati Andavari~~~
Manouchehr Zaker\footnote{Corresponding author: mzaker@iasbs.ac.ir}
\vspace{5mm}\\
   Department of Mathematics,\\
    Institute for Advanced Studies in Basic Sciences,\\
  Zanjan 45137-66731, Iran}
    \date{}
\begin{document}
\maketitle
\begin{abstract}
\noindent Dynamic monopolies were already defined and studied for the formulation of the phenomena
of the spread of influence in social networks such as disease, opinion, adaptation of new product and etc.
The elements of the network which have been influenced (e.g. infected or adapted an opinion) are called active
vertices. It is assumed in these models that when an element is activated, it remains active until the end of the process. But in some phenomena
of the spread of influence this property does not hold. For example in some diseases the infection lasts only a limited period
of time or consider the spread of disease or propagation of computer virus together with some quarantination or decontamination methods.
Dynamic monopolies are not useful for the study of these latter phenomena. For this purpose, we introduce
a new model for such diffusions of influence and call it weak dynamic monopoly. A social network is
represented by a graph $G$. Assume that any vertex $v$ of $G$ has a threshold $\tau(v)\in \Bbb{N}$.
Then a subset $D\subseteq V(G)$ is said to be a weak dynamic monopoly if $V(G)$ can be partitioned
into $D_0=D, D_1, \ldots, D_k$ such that for any $i$, any vertex $v$ of $D_i$ has at least $\tau(v)$
neighbors in $D_{i-1}$. In this definition, by the size and the processing time of $D$ we mean $|D|$ and $k$, respectively.
We first investigate the relationships between weak dynamic monopolies and
other related concepts and then obtain some bounds for the smallest size of weak dynamic monopolies. Next we obtain some results
concerning the processing time of weak dynamic monopolies in terms of some graph parameters. An upper bound is given for the smallest size of weak dynamic monopolies in the Cartesian product of cycles and its processing time is determined. Finally, a hardness result concerning inapproximibility of the determining the smallest size of weak dynamic monopolies in general graphs is obtained.
\end{abstract}




\section{Introduction}

\noindent The formulation and analysis of the various models of the spread of influence such as disease and opinion in a population, virus in computer networks, adaptation of innovation and viral marketing in social networks have been the research subject of many authors in the recent years.
Let a graph $G$ represent the underlying social network, where $V(G)$ denotes the individuals or the elements of the network and $E(G)$ represents the links or ties between them. These models and phenomena of the spread of influence have been studied using graph theory and in terms of (progressive) dynamic monopolies in graphs (or $k$-conversion process and target set selection in some other articles). Let a social network be represented by a graph $G=(V(G),E(G))$. We call such a graph $G$ social graph. By the social graph we only mean that the underlying graph corresponds to a social (or even virtual) network. We also assume that corresponding to any vertex
of $G$, there exists a threshold $\tau(v)$ such that $0\leq \tau(v)\leq deg(v)$, where
$deg(v)$ denotes the degree of $v$ in $G$. A subset $D\subseteq V(G)$ is said to be a dynamic monopoly
if $V(G)$ can be partitioned into $D_0=D, D_1, \ldots, D_k$ such that for any $i$, any vertex $v$ of
$D_i$ has at least $\tau(v)$ neighbors in $D_0\cup \ldots \cup D_{i-1}$. We say a vertex $v$ becomes active at time $i$ if
$v$ belongs to $D_i$. We denote by $\tau$ the threshold assignment for the vertices of $G$. The smallest cardinality of any dynamic monopoly of $(G,\tau)$ is denoted by $dyn_{\tau}(G)$. Dynamic monopolies are modeling the spread of influence in $G$, where $\tau(v)$ is interpreted as the degree of susceptibility of the vertex $v$. Dynamic monopolies were widely studied in the literature \cite{ABW,ABST,CL,C,DR,FKRRS,FLLPS,KSZ,Z,Z2}, under the equivalent term ``conversion sets" \cite{CDPRS, DR} and also ``target set selection" \cite{ABW, C, NNUW}. Dynamic monopolies have applications in viral marketing \cite{DRi}. In modeling the phenomena of the spread of influence, it is assumed that when an element is activated, it remains active forever. In practice some models of spread of influence such as spread of a disease together with quarantination or the spread of virus in computer networks combined with decontamination process, do not
match with the underlying conditions of dynamic monopolies \cite{F}. From the other side, it was known that in some diseases the infection lasts only a limited period of time for each infected person \cite{AM}.
Dynamic monopolies are not useful for the formulation of these latter phenomena. For this purpose, we introduce in this paper a new model for such diffusions of influence and call it weak dynamic monopoly. The following presents the formal definition and some related notations. For the graph theoretical terminology we refer the reader to \cite{BM}.
\begin{defin}
Assume that a graph $G$ together with a threshold assignment $\tau$ for its
vertices is given. A subset $D\subseteq V(G)$ is said to be a $\tau$-weak dynamic monopoly (or $\tau$-WDM) if $V(G)$ can
be partitioned into $D_0=D, D_1, \ldots, D_t$ such that for any $i$, any vertex $v$ of
$D_i$ has at least $\tau(v)$ neighbors in $D_{i-1}$. By the size of the weak dynamic
monopoly $D$ we mean the cardinality of $D$. Let $D$ be a weak dynamic monopoly in $(G,\tau)$ as above. We call the value $t$
in the above definition, the
processing time of $D$. We denote the smallest size of
any $\tau$-WDM of $(G,\tau)$ by $wdyn_{\tau}(G)$.
\end{defin}

\noindent A special case of weak dynamic monopolies, the so-called monopolies were already defined and studied
in the literature \cite{FFGHHKLS, KNSZ}. In a graph $G$, a subset $D\subseteq V(G)$ is said to be a monopoly
if each vertex $v \in V(G)\setminus D$ has at least $\tau(v)$ neighbors in $D$.
Denote the smallest size of any monopoly respect to the threshold function $\tau$ by $mo_{\tau}(G)$. Also Flocchini et al. in \cite{FKRRS}
study the dynamic monopolies which activate the whole vertices of the graph in only one time step and call them stamos. Note that stamos, monopolies,
and weak dynamic monopolies with the processing time one, are all equivalent concepts.

\noindent In the topics of monopolies or dynamic monopolies the following two special types of threshold assignments are mostly studied.
In {\it simple majority threshold} we have $\tau(v)=deg(v)/2$ for any vertex $v$ of the graph and in
{\it strict majority threshold}, the threshold of any vertex $v$ is defined as $\tau(v)=\lceil(deg(v)+1)/2\rceil$.

\noindent {\it The outline of the paper is as follows:} We end this section by comparing three quantities
$dyn(G)$, $mo(G)$ and $wdyn(G)$ (Theorem \ref{monopoly}). In Section 2 we obtain some bounds for the size of weak dynamic monopolies in terms of the order and the processing time or the even-girth of graphs for general threshold assignments. Section 3 devotes to study the processing time of weak dynamic monopolies in terms of known graph parameters. In Section 3, an upper bound is given for the smallest size of weak dynamic monopolies in the Cartesian product of cycles and its processing time is determined. Finally, in Section 4 we show that for any $\epsilon$, the smallest size of WDM can not be approximated within a factor of ${\mathcal{O}}(2^{log^{1-\epsilon}n})$, unless $NP\subseteq DTIME(n^{polylog(n)})$, where $n$ is the order of the input graph.

\noindent The following theorem gives the comparison results between monopoly, dynamic monopoly and weak dynamic monopoly. For any two non-negative valued functions $f(n)$ and $g(n)$, by $f(n)=o(g(n))$ we mean $f(n)/(g(n))$ tends to zero as $n\rightarrow \infty$. Also we write $f(n)=\Omega(g(n))$ if there exists a positive constant $\lambda$ such that for any $n$, $f(n)\geq \lambda g(n)$.

\begin{thm}\label{monopoly}

(i) Let $(G,\tau)$ be any graph. Then $dyn_{\tau}(G)\leq wdyn_{\tau}(G)\leq mo_{\tau}(G)$.

(ii) The equality may hold in each of the above inequalities.

(iii) There exists a sequence of graphs $H_n$ such that
$$dyn_{\tau}(H_n)=wdyn_{\tau}(H_n)=o(mo_{\tau}(H_n)).$$

(iv) There exists a sequence of graphs $G_n$ for which $$dyn_{\tau}(G_n) = o(wdyn_{\tau}(G_n)).$$
\end{thm}

\noindent\begin{proof} The validity of the inequalities in (i) is clear from the definitions. To prove (ii), it is enough to consider the
cycle graph $C_n$ on $n$ vertices with the strict majority threshold. We observe that
$dyn(C_n) = wdyn(C_n) = mo(C_n)=\lceil n/2 \rceil$. To prove (iii), replace any vertex of $C_n$ by $\overline{K_2}$ (the complement of the complete graph on two vertices) and then join
any two consecutive $\overline{K_2}$. Denote the resulting graph by $H_n$ and consider simple
majority threshold for $H_n$. It is easily seen that $dyn(H_n)=wdyn(H_n)=o(mo(H_n))$.

\noindent We now prove (iv). We obtain a sequence of graphs $G_n$ for which $dyn(G_n)=o(wdyn(G_n))$.
Define $G_n=K_1\vee C_n$, where $\vee$ is the join notation. Set
$\tau(v)=\lceil deg(v)/2 \rceil$ for any vertex $v$ of $G_n$. Denote the vertex of $K_1$
in $G_n$ by $u$. A set consisting of $u$ and one vertex from $C_n$ forms a dynamic monopoly.
Hence $dyn(G_n)=2$. We now show $wdyn(G_n)=\Omega(n)$. Let $D$ be any WDM of $G_n$. Assume first that $u$ does not
belong to $D$. In this case, since each vertex of $D_1$ has two neighbors in $D_0$, then $C_n\cap D_j=\emptyset$ for $j\geq 2$ and also $D_1$ is an independent set which implies
$|D|\geq \lceil n/2 \rceil$. Assume now that $u\in D$. Let $V(C_n)=\{0,1,2,\ldots,n-1\}$. In the following we show that the processing time
of $D$ is at most two. Suppose $D$ partitions $V(G_n)$ into $D_0=D, D_1, \ldots, D_t$,
where $t\geq 3$ and let $i\in D_3$. Denote the neighborhood of the vertex $i$ in $G_n$ by
$N(i)$. Since $\tau(i)=2$, we have necessarily $|N(i)\cap D_2|\geq 2$. And then $i+1, i-1\in D_2$.
Since $u\in D_0$ and $i\in D_3$, then $|N(i-1)\cap D_1|\leq 1$ and $i-1$ can not become active at time 2,
which is a contradiction.  Hence, the processing time is at most two. In the following we prove $|D_0|\geq n/4$.
Assume by the contrary that $|D_0|< n/4$. Then there must be four consecutive vertices $v_1, \ldots, v_4$ in $C_n$ which are not in $D_0$.
It is clear from the structure of $G_n$ that
no three consecutive vertices in $C_n$ are activated at time 1 and also $D_2$ is independent. Then there are three possibilities:

\noindent 1) $v_1,v_4\in D_2$ and $v_2,v_3\in D_1$. This case clearly rules out the activation of $v_2$ and $v_3$ at time 1, a contradiction.

\noindent 2) $v_1,v_3\in D_2$ and $v_2,v_4\in D_1$. In this case $v_2$ could not have two neighbors in $D_0$, a contradiction.

\noindent 3) $v_2\in D_2$ and $v_1,v_3,v_4\in D_1$. In this case $v_3$ could not have two neighbors in $D_0$, again a contradiction.

\noindent Therefore $|D_0|\geq n/4$. In other words, $wdyn(G_n)=\Omega(n)$, which completes the proof of (iv).
\end{proof}

\section{Some bounds for the size of weak dynamic monopolies}

\noindent We first consider weak dynamic monopolies with strict majority threshold.
\begin{thm}
For any graph $G$ with strict majority threshold $wdyn(G) \leq 2|G|/3$. Moreover the bound is tight.
\end{thm}

\noindent\begin{proof} It was shown in \cite{FFGHHKLS} that any graph $G$ contains a monopoly with strict majority threshold and no more than
$2|G|/3$ vertices. Then, by Theorem \ref{monopoly} we have $wdyn(G)\leq 2|G|/3$.

\noindent We now obtain a graph $G$ satisfying $wdyn(G)=2|G|/3$. Consider $m$ vertex disjoint copies of $K_3$. Add a new
vertex $v$ to the graph and then connect $v$ to exactly one vertex from each copy of $K_3$. Let $D$ be any minimum weak
dynamic monopoly in $G$, where the strict majority threshold is considered.
We claim that $D$ consists of two vertices from each triangle of $G$. Otherwise, let $v_1$, $v_2$ and $v_3$
be vertices of a triangle and $v_3$ be adjacent to the central vertex of the graph. One of the vertices
$v_1$ or $v_2$ has to be in $D$, otherwise, neither $v_1$ nor $v_2$ becomes active by the other vertices of
$G$. Assume that $v_1\notin D$. Since $v_1$ needs two active neighbors
at the same time, then necessarily $v_3\in D$. We have therefore
$$\lfloor\frac{2(3m+1)}{3}\rfloor=2m\leq wdyn(G)=|D|\leq \lfloor\frac{2(3m+1)}{3}\rfloor.$$
\end{proof}

\noindent It was proved in \cite{Z} that any strict majority dynamic monopoly in a cubic graph on $n$
vertices contains at least $(n+2)/4$ vertices. Hence by Theorem \ref{monopoly} we have the following.

\begin{thm}\label{lower} Let $G$ be a cubic graph on $n$ vertices and for each vertex $v$, $\tau(v)=2$. Then
$$\frac{n+2}{4}\leq wdyn(G).$$
\end{thm}
\noindent In Proposition \ref{tight} we show that the bound of Theorem \ref{lower} is tight. In the following we obtain a lower bound for the size of any weak dynamic monopoly in terms of
the processing time of the weak dynamic monopoly.

\begin{thm}\label{odd}
Let $G$ be a graph with the maximum degree $\Delta(G)\leq 2r+1$. Assume that each vertex $v$ of $G$ has threshold $\tau(v)$
such that $\tau(v)\geq r+1$. Suppose $D$ is any $\tau$-WDM for $G$ with processing time $t$. Then
$$\frac{n}{[2r+2-(2r+1)(\frac{r}{r+1})^t]}\leq |D|.$$
\end{thm}

\noindent\begin{proof} Let $D_0,\ldots ,D_t$ be a partition for $V(G)$ corresponding to WDM $D$ where $D_0=D$.
Let $e_i$, $i= 1 , \ldots , t$, be the number of edges between $D_{i-1}$ and $D_i$. Since
each vertex in $D_0$ has at most $2r+1$ neighbors in $D_1$ and each vertex in $D_1$ has
at least $r+1$ neighbors in $D_0$, then
$$(r+1)|D_1|\leq e_1 \leq (2r+1)|D_0|.$$

\noindent From the other side, for each $i\in \{1, \ldots ,t-1\}$, each vertex in $D_i$ has at most $r$
neighbors in $D_{i+1}$ and every vertex in $D_{i+1}$ has at least $r+1$ neighbors in $D_i$ . Therefore
$$(r+1)|D_{i+1}|\leq e_i \leq r|D_i|.$$
\noindent And hence
$$|D_{i+1}|\leq (\frac{r}{r+1})|D_i|.$$
\noindent Using the latter inequality repetitively, we obtain the following inequality
$$|D_i|\leq (\frac{2r+1}{r+1})(\frac{r}{r+1})^{i-1}|D_0|.$$
\noindent It turns out that
$$n=\sum_{i=0}^{t}|D_i|\leq |D_0|\left(1+\sum_{i=1}^{t}(\frac{2r+1}{r+1})(\frac{r}{r+1})^{i-1}\right)$$
\noindent which implies
$$\frac{n}{[2r+2-(2r+1)(\frac{r}{r+1})^t]}\leq |D_0|.$$
\end{proof}

\noindent In order to present Proposition \ref{tight} we need the following easy number theoretic fact.

\begin{lemma}\label{number}
For any positive odd integer $t\geq 3$, $2^{t-1}-1$ is divisible by $3$.
\end{lemma}

\noindent \begin{proof}
Let $t=2p+1$ for some integer $p\geq 1$. The proof is easily obtained by the induction on $p$.
\end{proof}

\noindent The following proposition shows that in Theorem \ref{lower} and Theorem \ref{odd} the equality may hold for some graphs.

\begin{prop}\label{tight}
For any integers $k$ and $t$ with $3k+1=2^{t-1}$, there exists a cubic graph $G$ on $n=8k+2$ vertices which
contains a weak dynamic monopoly $D$ corresponding to the constant threshold assignment $2$ for the vertices of $G$, such that the processing time of $D$ is $t$ and $|D|= (n+2)/4$.
\end{prop}

\noindent\begin{proof} By Lemma \ref{number} there are infinitely many $t$ for which there exists $k$ such that $3k+1=2^{t-1}$. Let $k$ and $t$ be two arbitrary integers such that $3k+1=2^{t-1}$. We construct a
cubic graph $G$ on $n=8k+2$ vertices in which each vertex has threshold 2 and
$$wdyn(G)=\left\lceil \frac{n}{\left( 2r+2-(2r+1)(r/(r+1))^t \right)}\right\rceil=\frac{n+2}{4}.$$

\noindent Note that in the above relation $r=1$, since the graph is cubic. It can be easily shown that the second equality in the above relation holds for the values $n=(2^{t+2}-2)/3$ and $r=1$. Now we explain the construction of $G$. Consider a partition for $V(G)$ in the form of $V(G)=D_0\cup D_1\cup \ldots\cup D_t$ such that $|D_0|=2k+1$ and
$|D_1|=3k+1=2^{t-1}$ and for each $i$, $i\in\{2,\ldots,t\}$, $|D_i|=2^{t-i}$. Note that $|V(G)|=2^t-1+|D_0|=8k+2$. Put $3(2k+1)-1$ edges between
$D_0$ and $D_1$ in such a way that each vertex in $D_1$ has exactly two neighbors in $D_0$ and each vertex in
$D_0$ except one vertex say $v$, has three neighbors in $D_1$. The vertex $v$ of $D_0$ has only two neighbors in $D_1$.
For each $i$, $i\geq 1$, put some edges between $D_i$ and $D_{i+1}$ in such a way that each vertex in $D_i$ has exactly one neighbor in $D_{i+1}$ and each vertex in $D_{i+1}$
has exactly two neighbors in $D_i$. Now connect the only vertex of $D_t$ to the two vertices in $D_{t-1}$ and also to
$v\in D_0$. Consider the constant threshold assignment $2$ for all the vertices of $G$. It is easily seen that $D_0$ is a WDM with the desired cardinality. Therefore in the inequalities of Theorem \ref{lower} and Theorem \ref{odd}, equality hold for the graph $G$.
\end{proof}

\noindent For $2r$-regular graphs we have the following analogous bound.

\begin{thm}\label{even}

\noindent (i) Let $G$ be a $2r$-regular graph on $n$ vertices and set $\tau(v)=r+1$ for every vertex $v$ of $G$.
Let also $D$ be any $\tau$-WDM with processing time $t$. Then
$$\frac{n}{[1+r(1-(\frac{r-1}{r+1})^t)]}\leq |D|.$$
(ii) If $\tau(v)=r$ for each vertex $v$, then
$$\frac{n}{1+2t}\leq |D|.$$
\end{thm}

\noindent\begin{proof}

\noindent (i) If $r=1$ then the graph is the cycle $C_n$. Using Theorem \ref{monopoly} we have
$$wdyn(C_n)=dyn(C_n)=mo(C_n)=\lceil n/2 \rceil\geq \frac{n}{[1+1(1-(\frac{1-1}{1+1})^1)]}.$$
For $r\geq 2$, the proof is similar to the proof of Theorem \ref{odd}.

\noindent (ii) For this case too, the proof is similar to the proof of Theorem \ref{odd}.
\end{proof}


\noindent The following theorem presents a lower bound in terms of even-girth of the graph. By the even-girth, we mean the length of smallest even cycle in the graph. In the following by $d(x,y)$, for any two vertices $x$ and $y$, we mean the distance between $x$ and $y$ in the graph.

\begin{thm}\label{eg}
Let $G$ be a graph with even-girth $eg(G)=2k+2$, where $k\geq 3$. Let $\tau$ be a threshold assignment
for the vertices of $G$ such that for any vertex $v$, $0 < \tau(v) < deg(v)$ and $t_m= {\min}_{v\in G} \tau(v)\geq 3$. Then
$$(t_m-1)^{\lfloor\frac{k}{2}\rfloor+1}\leq wdyn_{\tau}(G).$$
\end{thm}

\noindent\begin{proof} Let $n$ be the order of $G$. To prove the theorem, we need to show the following inequality.

$$n=|G|\geq 1+\frac{(t_m+1)}{(t_m-2)}\big((t_m-1)^k-1)\big)+t_m(t_m-2).$$

\noindent For this purpose, consider some vertex $v\in V(G)$ and let $N_i=\{u\in V(G):d(u,v)=i\}$. Since
$eg(G)=2k+2$ then each vertex $u\in N_i$, $1\leq i\leq k$, has at most one neighbor in $N_i$ and
at most one neighbor in $\{v\}\cup (\cup_{j=1}^{i-1}N_i)$.
Hence each $u\in N_i$, $1\leq i\leq k-1$, has at least $t_m-1$ neighbors in $N_{i+1}$.
We have also $deg(v)\geq t_m+1$ for each vertex $v\in V(G)$. Therefore $|N_1|\geq (t_m+1)$ and then
$|N_i|\geq (t_m+1)(t_m-1)^{i-1}$.

\noindent Now suppose $u\in N_{k-1}$ and $\{u_1,u_2,\ldots,u_{t_m-1}\}\subseteq N(u)\cap N_k$, where $N(u)$ is the neighborhood of $u$ in $G$.
Since $eg(G)=2k+2\geq 8$ then for any $1\leq i,j\leq t_m-1$, $u_i$ and $u_j$ have no common
neighbor. Let $\{w_1,\ldots,w_{t_m+1}\}\subseteq N_1$ and $C_{w_j}=\{z\in N_k : d(w_j,z)=k-1\}$, where
$1\leq j\leq t_m+1$. Let also $\{u_1,u_2,\ldots,u_{t_m-1}\}\subseteq C_{w_1}$.
If $u_1$ is adjacent to some vertex in $C_{w_j}$ with $j\neq 1$, then none
of $u_2,\ldots,u_{t_m-1}$ have any neighbor in $C_{w_j}$. Now suppose without loss
of generality that $u_1$ is adjacent to some vertex from each $C_{w_j}$, $2\leq j\leq t_m+1$. Since for
each $i\in \{2,\ldots,t_m-1\}$ we have $deg(u_i)\geq t_m+1$, then each of these vertices has at least $t_m$ neighbors out of $\{v\}\cup (\cup_{i=1}^{k}N_i)$.
Therefore we have
$$n\geq 1+\big(\sum_{i=1}^{k}(t_m+1)(t_m-1)^{i-1}\big)+t_m(t_m-2).$$

\noindent Now let $D_0$ be any WDM with the processing time $t$. The set $D_0$ partitions $V(G)$ as
$V(G)=\cup_{i=0}^t D_i$. There are two possibilities for $t$:

\noindent{\bf 1)} $t\geq \lceil k/2\rceil+1$.

\noindent In this case we consider a vertex say $v\in D_t$. The vertex $v$ has at least $t_m$ neighbors in $D_{t-1}$.
Each of these neighbors has at least $t_m$ neighbors in $D_{t-2}$ and in general for any $0<i\leq t$ any vertex in $D_i$ has at least $t_m$ neighbors in $D_{i-1}$. Then
$$|D_0|\geq t_m^t\geq t_m^{\lceil k/2\rceil+1}\geq t_m^{\lfloor k/2\rfloor+1}\geq (t_m-1)^{\lfloor k/2\rfloor+1}.$$

\noindent{\bf 2)} $t\leq \lceil k/2\rceil$.

\noindent In this case we have the following two subcases.

\noindent{\bf Subcase \textit{i})} $\forall~ i\in\{0,\ldots,t-1\} ; |D_i|\geq |D_{i+1}|$.

\noindent In this subcase we have
$$n=\sum_{i=0}^t |D_i|\leq (t+1)|D_0|.$$
Hence
\begin{align*}
|D_0| & \geq\frac{n}{t+1}\\
    &\geq\frac{n}{\lceil\frac{k}{2}\rceil+1}\\
    &\geq \frac{1}{\lceil\frac{k}{2}\rceil+1}\bigg(1+\frac{(t_m+1)}{(t_m-2)}\big((t_m-1)^k-1\big)+t_m(t_m-2)\bigg).
\end{align*}

\noindent For $k\in \{3,4\}$, it is easily seen that $|D_0|\geq (t_m-1)^{\lfloor k/2\rfloor+1}$. When $k\geq 5$ we have

\begin{align*}
|D_0| &\geq \frac{1}{\lceil\frac{k}{2}\rceil+1}\bigg(1+\frac{(t_m+1)}{(t_m-2)}\big((t_m-1)^k-1\big)+t_m(t_m-2)\bigg)\\
    &\geq \frac{1}{\lceil\frac{k}{2}\rceil+1}\bigg(1+(t_m-1)^k-\frac{(t_m+1)}{(t_m-2)}+t_m(t_m-2)\bigg)\\
    &\geq\frac{(t_m-1)^k}{\lceil\frac{k}{2}\rceil+1}\\
    &=(t_m-1)^{\lfloor\frac{k}{2}\rfloor+1}\bigg(\frac{(t_m-1)^{\lceil\frac{k}{2}\rceil-1}}{\lceil\frac{k}{2}\rceil+1}\bigg)\\
    &\geq(t_m-1)^{\lfloor\frac{k}{2}\rfloor+1}\bigg(\frac{2^{\lceil\frac{k}{2}\rceil-1}}{\lceil\frac{k}{2}\rceil+1}\bigg).
\end{align*}

\noindent It is easily seen that for $k\geq 5$, $2^{\lceil k/2\rceil-1}/(\lceil k/2\rceil+1)\geq 1$. Then
$$|D_0|\geq (t_m-1)^{\lfloor k/2\rfloor+1}.$$

\noindent{\bf Subcase \textit{ii})} $\exists ~ i\in\{0,\ldots,t-1\}; |D_i| < |D_{i+1}|$.

\noindent Let $i$ be the smallest one with this property. Consider the bipartite graph constructed on
sets $D_i$ and $D_{i+1}$ and all edges of $G$ between them. Denote this new graph by $H$. Then $|D_{i+1}|>|N_H(D_{i+1})|$. Suppose
$S$ is a minimal subset of $D_{i+1}$ with the property $|S|>|N_H(S)|$. It is easily seen that for each vertex $v\in N_H(S)$, $|N(v)\cap S|=deg_S(v)\geq 2$.
Otherwise, let $v\in N_H(S)$ has only one neighbor in $S$. Set $S'=S\setminus N_H(v)$. Then
$|S'|=|S\setminus N_H(v)|=|S|-1>|N_H(S)|-1=|N_H(S')|$, a contradiction with the minimality of $S$.
On the other hand, every vertex in $S$ has at least $t_m$ neighbors in $N_H(S)$. Then we have
$$\sum_{v\in S}deg_{N_H(S)}(v)=\sum_{u\in N_H(S)}deg_S(u)\geq t_m|S|\geq t_m(|N_H(S)|+1).$$
Hence, there is a vertex say $u$ in $N_H(S)$, such that $deg_S(u)\geq t_m+1$. Consider induced subgraph of $H$ on $S$ and $N_H(S)$ as
$H^S=H[S,N_H(S)]$ and  $v\in N_S(u)$. Let $N_i^S=\{w\in V(H^S) : d_{H^S}(v,w)=i\}$. Now we have $|N^S_1|\geq t_m$. Since $deg_S(u)\geq t_m+1$
and every vertex in $N_H(S)$ has at least two neighbors in $S$ then $|N_2^S|\geq t_m+t_m-1=2t_m-1$. It is easily seen that $N_2^S\subseteq S$
and each vertex in $S$ has at least $t_m$ neighbor in $N_H(S)$. Then
$|N_3^s|\geq (2t_m-1)(t_m-1)$. By repeating this argument we obtain that for even $k$, $N_k^S\subseteq S$ and $|N_k^S|\geq (2t_m-1)(t_m-1)^{k/2-1}$. Since the even girth of $G$ is $2k+2$, then the girth of $H$ is at least $2k+2$. It implies that any vertex of $N^S_k$ has at most
one neighbor in $\cup_{i<k}N_i^S$. Since $H$ is bipartite then there must be $(2t_m-1)(t_m-1)^{k/2}$ neighbors for the vertices of $N^S_k$ in $N_H(S)\setminus \cup_{i<k}N_i$. Therefore
$$|D_i|\geq |N_H(S)|\geq (2t_m-1)(t_m-1)^{\frac{k}{2}}\geq (t_m-1)^{\frac{k}{2}+1}=(t_m-1)^{\lfloor\frac{k}{2}\rfloor+1}.$$
For odd $k$, $N_k^S\subseteq N_H(S)$ and $|N_k^S|\geq (2t_m-1)(t_m-1)^{\lfloor\frac{k}{2}\rfloor}$. Then we have
\begin{align*}
|N_H(S)|&\geq t_m+(2t_m-1)(t_m-1)+\ldots+(2t_m-1)(t_m-1)^{\lfloor\frac{k}{2}\rfloor}\\
        &= t_m+(2t_m-1)\big(\frac{(t_m-1)^{\lfloor\frac{k}{2}\rfloor+1}-1}{t_m-2}-1\big)\\
        &\geq t_m+2(t_m-1)^{\lfloor\frac{k}{2}\rfloor+1}-2-2t_m+4\\
        &\geq (t_m-1)^{\lfloor\frac{k}{2}\rfloor+1}.
\end{align*}
Hence $|D_0|\geq|D_i|\geq |N_H(S)|\geq (t_m-1)^{\lfloor\frac{k}{2}\rfloor+1}$. This completes the proof.
\end{proof}

\section{Bounds for processing time}

\noindent In the following we consider the processing time of any WDM and obtain some upper bounds for it.

\begin{thm}\label{paths} Let $G$ be a graph with a threshold assignment $\tau$, where
$\tau(v)\geq k\geq 1$ for each vertex $v$ and let $D_0$ be any WDM which partitions $V(G)$ as $\cup_{i=0}^tD_i$. Then
there are $k$ internally disjoint paths with length $t-1$  beginning from
$D_0$ and containing only one vertex from each $D_i$, $0\leq i \leq t-1$.
\end{thm}

\noindent\begin{proof} By induction on $i$, where $0\leq i\leq t-1$, we show that for any set
$\{u_1,\ldots,u_k\}\subseteq D_i$ of distinct vertices,
there are $k$ internally disjoint paths with length $i$ which begin from $D_0$ and end at
$u_1,\ldots,u_k$, respectively and these paths  contain only one vertex from each $D_j$, $0\leq j < i$.
This is obvious when $\{u_1,\ldots,u_k\}\subseteq D_0$. Suppose $\{u_1,\ldots,u_k\}\subseteq D_i$, $0 < i \leq t-1$.
Because each of $u_1,\ldots,u_k$ has at least $k$ neighbors in $D_{i-1}$, then there are distinct
$u_1',\ldots,u'_k$ in $D_{i-1}$ such that $u'_s \in N(u_s)$, $1\leq s \leq k$. By the induction
hypothesis  there are internally disjoint  paths $P'_1,\ldots,P'_k$ with length $i-1$ which begin
from $D_0$ and end at $u'_1,\ldots,u'_k$, respectively and contain only one vertex from each $D_j$,
$j < i$. Now, $P'_1+u'_1u_1,\ldots,P'_k+u'_ku_k$ are the desired paths, where by $P'_1+u'_1u_1$ we mean the extension of path $P'_1$ by the edge $u'_1u_1$.
\end{proof}

\noindent The following corollaries are obtained from Theorem \ref{paths}. By a starlike tree we mean
any tree that is isomorphic to a subdivision of $K_{1,n}$ for some $n$. Such a starlike tree contains a central vertex of degree $n$ and $n$ branches.

\begin{cor}\label{starlike}
Let $G$ be a graph with threshold assignment $\tau$, where for each vertex $v$ of $G$,
$\tau(v)\geq k\geq 2$. Then there exists a starlike tree in $G$  with
central vertex of degree $k$ and $k$ branches of length $t-1$.
\end{cor}

\noindent\begin{proof} Let $v\in D_t$. The vertex $v$ has at least $k$ neighbors $u_1,\ldots,u_k$ in $D_{t-1}$.
By Theorem \ref{paths}, there are $k$ distinct paths $P_1,\ldots,P_k$ of length $t-1$
that end at $u_1,\ldots,u_k$, respectively. A starlike
tree can be easily obtained by adding vertex $v$ and edges $vu_1,\ldots,vu_k$ to $P_1,\ldots,P_k$.
\end{proof}

\noindent In the following corollary by $\alpha'(G)$ we mean the maximum number of independent edges in $G$.

\begin{cor}
Let $G$ be a graph with threshold assignment $\tau$ and $\min_{v\in V(G)} \tau(v)\geq k$. Then for any WDM with
processing time $t$ we have
$$t\leq\frac{2\alpha'(G)}{k}+2.$$
\end{cor}

\noindent\begin{proof} Since $\min \tau(v)\geq k$, by Corollary \ref{starlike} there is a
starlike tree with $k$ branches of length $t-1$. By choosing $\lfloor(t-1)/2\rfloor$ independent edges
from each branch we have a set containing $k(t-1)/2$ independent edges. Therefore
$$k(\frac{t-2}{2})\leq k\lfloor\frac{t-1}{2}\rfloor\leq \alpha'(G).$$
This yields the desired result.
\end{proof}

\noindent In the following corollary we present an upper bound for processing time in terms of the length of longest path in graphs.

\begin{cor}
Let $G$ be a graph with threshold assignment $\tau$ and let $l$ be the length of longest path in $G$. Then
for any WDM with processing time $t$ we have

\item{ i) If for each vertex $v$, $\tau(v)\geq 2$   then  $t\leq l/2$};
\item{ ii) If $\tau$ is the strict majority assignment, then  $t\leq (l+ 2)/2$}.
\end{cor}

\noindent\begin{proof} i) By corollary \ref{starlike}, there is a path of length $2(t-1)+2=2t$ in $G$. Hence $l\geq 2t$.

\noindent ii) Let $D$ be a WDM which partitions $V(G)$ as  $\cup_{i=1}^tD_i$. Let also $w\in D_t$.
Then $w$ has at least one neighbor $w'$ in $D_{t-1}$ with $deg(w')\geq 2$.
Since the activation process follows the strict majority rule, then $w'$ has at least two neighbors in $D_{t-2}$ such as
$v$ and $u$ whose degrees and thresholds in $G$ are at least two. By Theorem \ref{paths} there
exist two internally disjoint paths $P_1$ and $P_2$ with length $t-2$ which end at $v$ and $u$,
respectively. Now $P_1+vw'+P_2+uw'$ is a path of length $2t-2$. Then $t\leq(l+2)/2$.
\end{proof}


\noindent We are going to show that in graphs $G$ with bounded maximum degree, the size of any WDM for $G$ and its processing
time are not bounded by a constant value, i.e. one of them goes to infinity as $|G|\rightarrow \infty$.

\begin{thm}\label{tdelta}
Let $G$ be a graph on $n$ vertices such that $\Delta(G)\leq k$ for some constant $k$.
Let $\tau$ be any threshold assignment for $G$ and $D$ be any $\tau$-WDM with the processing time $t$.
Then
$$n\leq k^{t+1}|D|.$$
\end{thm}

\noindent\begin{proof} Assume that $D$ partitions $V(G)$ as $D\cup D_1\cup\ldots\cup D_t$. It is easily seen that for
each $i$, $i\in\{1,\ldots,t\}$, $|D_i|\leq\Delta(G)^i |D_0|$. Then
$$n=\sum_{i=0}^t|D_i|\leq\sum_{i=0}^t \Delta(G)^i |D_0|\leq \sum_{i=0}^tk^i|D|.$$
Therefore
$$n\leq k^{t+1}|D|.$$
\end{proof}


\noindent It can be shown that Theorem \ref{tdelta} is not valid when $\Delta(G)$ is not bounded. For example consider $K_{2k}\vee C_n$ with any
threshold assignment $\tau$ such that $k\leq {\max}_{v\in V(C_n)}\tau(v)\leq 2k$ and ${\max}_{v\in V(K_{2k})}\tau(v)\leq n$. It is easily seen that any minimum WDM in this graph is
contained in $V(K_{2k})$. Hence $wdyn(K_{2k}\vee C_n)\leq 2k$. From the other side, in such a minimum WDM, all the vertices are activated in at most two steps.


\noindent In the following theorem we determine the processing time in the Cartesian product of cycles denoted by $C_n\Box C_m$. And also we obtain an upper bound for the smallest size of weak dynamic monopolies. Dynamic monopolies of this family of graphs were studied in \cite{FLLPS}.

\begin{thm} Let $G=C_n\Box C_m$, where the threshold of each vertex is 3.
Then, the activation process for any WDM in $G$ ends after two steps and $wdyn(C_n\Box C_n)\leq 3n^2/8$
when $4|n$.
\end{thm}

\noindent\begin{proof} Let $D$ be any WDM in $C_n\Box C_m$ and assume on the contrary that a vertex say $v$
becomes active at time 3. Then $D$ partitions $V(C_n\Box C_m)$ as $D=D_0, D_1,\ldots,D_t$ for some $t$,
$t\geq 3$. The vertex $v$ has three neighbors $v_1$, $v_2$ and $v_3$ in
$D_2$. Let $u\in D_1$ be a vertex which is adjacent to two neighbors of $v$. The vertex $u$ should have three neighbors in $D_0$
but this is impossible. Hence any WDM of $G$ activates the whole graph in at most two time steps.

\noindent In Figure \ref{indepWDM}, a WDM for $C_n\Box C_n$ which is also an independent set is presented. Let $D_0$ be an independent WDM in $C_n\Box C_n$. Let also the WDM $D_0$, partitions the graph into three sets $D_0, D_1, D_2$.
Obviously, each vertex in $D_1$ has at least three neighbors in $D_0$, then we have
$$3|D_1|\leq 4|D_0|.$$
\noindent We make the following claims concerning $D_1$ and $D_2$.

\noindent{\bf Claim 1)} $D_1$ is independent.

\noindent Otherwise, we obtain Figure 1, in which there are
adjacent vertices in $D_0$. This is a contradiction.

\vspace{0.3cm}
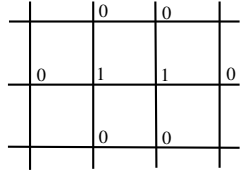
\begin{figure}[h]
\centering
\scalebox{0.7} 
{
\begin{pspicture}(0,-1.62)(4.42,1.62)
\psline[linewidth=0.04cm](0.42,1.6)(0.42,-1.6)
\psline[linewidth=0.04cm](1.62,1.6)(1.62,-1.58)
\psline[linewidth=0.04cm](2.82,-1.56)(2.8,1.6)
\psline[linewidth=0.04cm](4.02,1.6)(4.02,-1.56)
\psline[linewidth=0.04cm](4.4,1.22)(0.02,1.22)
\psline[linewidth=0.04cm](0.0,0.02)(4.4,0.02)
\psline[linewidth=0.04cm](4.4,-1.18)(0.0,-1.16)
\usefont{T1}{ptm}{m}{n}
\rput(1.766875,0.23){1}
\usefont{T1}{ptm}{m}{n}
\rput(2.986875,0.21){1}
\usefont{T1}{ptm}{m}{n}
\rput(1.8170313,-0.97){0}
\usefont{T1}{ptm}{m}{n}
\rput(3.0170312,-0.99){0}
\usefont{T1}{ptm}{m}{n}
\rput(4.2370315,0.23){0}
\usefont{T1}{ptm}{m}{n}
\rput(3.0170312,1.39){0}
\usefont{T1}{ptm}{m}{n}
\rput(1.8170313,1.43){0}
\usefont{T1}{ptm}{m}{n}
\rput(0.63703126,0.21){0}
\end{pspicture}
}
\caption{Two adjacent vertices in $D_1$.}
\end{figure}

\noindent{\bf Claim 2)} $D_2$ is independent and no vertex in $D_2$ has a neighbor in $D_0$.

\noindent If we have
two adjacent vertices in $D_2$ then the specified vertices in part (1) of Figure 2 can not
have three neighbors in the previous time.
And also we have $D_2\cap N(D_0)=\emptyset$. Otherwise, we have part (2) of Figure 2 which contradicts the independency of $D_0$.

\vspace{0.3cm}
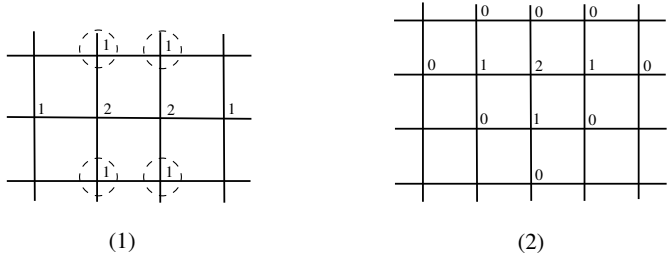
\begin{figure}[h]
\centering
\scalebox{0.6} 
{
\begin{pspicture}(0,-2.8314064)(14.64,2.8114061)
\usefont{T1}{ptm}{m}{n}
\rput(11.797031,-1.0385938){0}
\usefont{T1}{ptm}{m}{n}
\rput(10.577031,0.20140626){0}
\usefont{T1}{ptm}{m}{n}
\rput(9.417031,1.4014063){0}
\usefont{T1}{ptm}{m}{n}
\rput(11.746875,0.22140625){1}
\usefont{T1}{ptm}{m}{n}
\rput(10.566875,1.4014063){1}
\usefont{T1}{ptm}{m}{n}
\rput(10.597032,2.6014063){0}
\usefont{T1}{ptm}{m}{n}
\rput(11.7985935,1.3814063){2}
\usefont{T1}{ptm}{m}{n}
\rput(12.977032,0.18140624){0}
\usefont{T1}{ptm}{m}{n}
\rput(12.966875,1.4014063){1}
\usefont{T1}{ptm}{m}{n}
\rput(14.197031,1.3814063){0}
\usefont{T1}{ptm}{m}{n}
\rput(12.997031,2.5814064){0}
\usefont{T1}{ptm}{m}{n}
\rput(11.797031,2.5814064){0}
\psline[linewidth=0.04cm](8.6,-1.2285937)(14.62,-1.2285937)
\psline[linewidth=0.04cm](8.6,-0.00859375)(14.62,-0.00859375)
\psline[linewidth=0.04cm](8.58,1.1914062)(14.6,1.1914062)
\psline[linewidth=0.04cm](8.58,2.3914063)(14.6,2.3914063)
\psline[linewidth=0.04cm](9.22,2.7914062)(9.22,-1.6285938)
\psline[linewidth=0.04cm](10.4,2.7514062)(10.42,-1.6285938)
\psline[linewidth=0.04cm](11.6,2.7314062)(11.62,-1.6085937)
\psline[linewidth=0.04cm](12.8,2.7514062)(12.8,-1.6085937)
\psline[linewidth=0.04cm](14.0,2.7514062)(14.0,-1.5685937)
\psline[linewidth=0.04cm](5.42,1.6114062)(0.0,1.6114062)
\psline[linewidth=0.04cm](5.42,0.21140625)(0.0,0.25140625)
\psline[linewidth=0.04cm](0.0,-1.1685938)(5.4,-1.1685938)
\psline[linewidth=0.04cm](4.8,2.0914063)(4.82,-1.6485938)
\psline[linewidth=0.04cm](3.4,2.1314063)(3.4,-1.6485938)
\psline[linewidth=0.04cm](2.0,2.1114063)(2.0,-1.6285938)
\psline[linewidth=0.04cm](0.6,2.1514063)(0.62,-1.6085937)
\usefont{T1}{ptm}{m}{n}
\rput(2.2385938,0.44140625){2}
\usefont{T1}{ptm}{m}{n}
\rput(3.6385937,0.42140624){2}
\usefont{T1}{ptm}{m}{n}
\rput(3.586875,1.8414062){1}
\usefont{T1}{ptm}{m}{n}
\rput(2.206875,1.8614062){1}
\usefont{T1}{ptm}{m}{n}
\rput(0.766875,0.44140625){1}
\usefont{T1}{ptm}{m}{n}
\rput(2.206875,-0.9585937){1}
\usefont{T1}{ptm}{m}{n}
\rput(3.606875,-0.9585937){1}
\usefont{T1}{ptm}{m}{n}
\rput(4.986875,0.44140625){1}
\pscircle[linewidth=0.024,linestyle=dashed,dash=0.16cm 0.16cm,dimen=outer](3.45,1.7414062){0.43}
\pscircle[linewidth=0.024,linestyle=dashed,dash=0.16cm 0.16cm,dimen=outer](2.03,1.7614063){0.43}
\pscircle[linewidth=0.024,linestyle=dashed,dash=0.16cm 0.16cm,dimen=outer](3.45,-1.0785937){0.43}
\pscircle[linewidth=0.024,linestyle=dashed,dash=0.16cm 0.16cm,dimen=outer](2.03,-1.0785937){0.43}
\usefont{T1}{ptm}{m}{n}
\rput(2.5603125,-2.5285938){\Large (1)}
\usefont{T1}{ptm}{m}{n}
\rput(11.620313,-2.5485938){\Large (2)}
\end{pspicture}
}
\caption{(1) Two adjacent vertices in $D_2$. (2) Some vertex in $D_2$ is adjacent to a vertex in $D_0$.}
\end{figure}

\noindent Hence, from the above facts we have $|D_2|\leq(1/4)|D_1|$ and then $|D_2|\leq (1/3)|D_0|$. Therefore
$$mn=\sum_{i=0}^2 |D_i|\leq (1+\frac{1}{3}+\frac{4}{3})|D_0|.$$
\noindent And finally $|D_0|\geq 3nm/8$.

\begin{figure}\label{indepWDM}
\centering
\scalebox{0.6} 
{
\begin{pspicture}(0,-5.62)(11.62,5.62)
\psline[linewidth=0.04cm](0.8,4.2)(10.4,4.2)
\psline[linewidth=0.04cm](0.8,3.0)(10.4,3.0)
\psline[linewidth=0.04cm](0.8,1.8)(10.4,1.8)
\psline[linewidth=0.04cm](0.8,0.6)(10.4,0.6)
\psline[linewidth=0.04cm](0.8,-0.6)(10.4,-0.6)
\psline[linewidth=0.04cm](0.8,-1.8)(10.4,-1.8)
\psline[linewidth=0.04cm](0.8,-3.0)(10.4,-3.0)
\psline[linewidth=0.04cm](0.8,-4.2)(10.4,-4.2)
\psline[linewidth=0.04cm](9.8,4.6)(9.8,-4.6)
\psline[linewidth=0.04cm](8.6,4.6)(8.6,-4.6)
\psline[linewidth=0.04cm](7.4,4.6)(7.4,-4.6)
\psline[linewidth=0.04cm](6.2,4.6)(6.2,-4.6)
\psline[linewidth=0.04cm](5.0,4.6)(5.0,-4.6)
\psline[linewidth=0.04cm](3.8,4.6)(3.8,-4.6)
\psline[linewidth=0.04cm](2.6,4.6)(2.6,-4.6)
\psline[linewidth=0.04cm](1.4,4.6)(1.4,-4.6)
\usefont{T1}{ptm}{m}{n}
\rput(1.646875,4.51){1}
\usefont{T1}{ptm}{m}{n}
\rput(1.6770313,3.31){0}
\usefont{T1}{ptm}{m}{n}
\rput(2.8770313,4.51){0}
\usefont{T1}{ptm}{m}{n}
\rput(1.646875,2.11){1}
\usefont{T1}{ptm}{m}{n}
\rput(2.846875,3.31){1}
\usefont{T1}{ptm}{m}{n}
\rput(4.046875,4.51){1}
\usefont{T1}{ptm}{m}{n}
\rput(1.6770313,0.91){0}
\usefont{T1}{ptm}{m}{n}
\rput(2.8785937,2.11){2}
\usefont{T1}{ptm}{m}{n}
\rput(4.077031,3.31){0}
\usefont{T1}{ptm}{m}{n}
\rput(5.2785935,4.51){2}
\usefont{T1}{ptm}{m}{n}
\rput(1.646875,-0.29){1}
\usefont{T1}{ptm}{m}{n}
\rput(2.846875,0.91){1}
\usefont{T1}{ptm}{m}{n}
\rput(4.046875,2.11){1}
\usefont{T1}{ptm}{m}{n}
\rput(5.246875,3.31){1}
\usefont{T1}{ptm}{m}{n}
\rput(6.446875,4.51){1}
\usefont{T1}{ptm}{m}{n}
\rput(1.6770313,-1.49){0}
\usefont{T1}{ptm}{m}{n}
\rput(2.8770313,-0.29){0}
\usefont{T1}{ptm}{m}{n}
\rput(4.077031,0.91){0}
\usefont{T1}{ptm}{m}{n}
\rput(5.2770314,2.11){0}
\usefont{T1}{ptm}{m}{n}
\rput(6.477031,3.31){0}
\usefont{T1}{ptm}{m}{n}
\rput(7.677031,4.51){0}
\usefont{T1}{ptm}{m}{n}
\rput(1.646875,-2.69){1}
\usefont{T1}{ptm}{m}{n}
\rput(2.846875,-1.49){1}
\usefont{T1}{ptm}{m}{n}
\rput(4.046875,-0.29){1}
\usefont{T1}{ptm}{m}{n}
\rput(5.246875,0.91){1}
\usefont{T1}{ptm}{m}{n}
\rput(6.446875,2.11){1}
\usefont{T1}{ptm}{m}{n}
\rput(7.646875,3.31){1}
\usefont{T1}{ptm}{m}{n}
\rput(8.846875,4.51){1}
\usefont{T1}{ptm}{m}{n}
\rput(1.6770313,-3.89){0}
\usefont{T1}{ptm}{m}{n}
\rput(2.8785937,-2.69){2}
\usefont{T1}{ptm}{m}{n}
\rput(4.077031,-1.49){0}
\usefont{T1}{ptm}{m}{n}
\rput(5.2785935,-0.29){2}
\usefont{T1}{ptm}{m}{n}
\rput(6.477031,0.91){0}
\usefont{T1}{ptm}{m}{n}
\rput(7.6785936,2.11){2}
\usefont{T1}{ptm}{m}{n}
\rput(8.877031,3.31){0}
\usefont{T1}{ptm}{m}{n}
\rput(10.078594,4.51){2}
\usefont{T1}{ptm}{m}{n}
\rput(2.846875,-3.89){1}
\usefont{T1}{ptm}{m}{n}
\rput(4.046875,-2.69){1}
\usefont{T1}{ptm}{m}{n}
\rput(5.246875,-1.49){1}
\usefont{T1}{ptm}{m}{n}
\rput(6.446875,-0.29){1}
\usefont{T1}{ptm}{m}{n}
\rput(7.646875,0.91){1}
\usefont{T1}{ptm}{m}{n}
\rput(8.846875,2.11){1}
\usefont{T1}{ptm}{m}{n}
\rput(10.046875,3.31){1}
\usefont{T1}{ptm}{m}{n}
\rput(4.077031,-3.89){0}
\usefont{T1}{ptm}{m}{n}
\rput(5.2770314,-2.69){0}
\usefont{T1}{ptm}{m}{n}
\rput(6.477031,-1.49){0}
\usefont{T1}{ptm}{m}{n}
\rput(7.677031,-0.29){0}
\usefont{T1}{ptm}{m}{n}
\rput(8.877031,0.91){0}
\usefont{T1}{ptm}{m}{n}
\rput(10.077031,2.11){0}
\usefont{T1}{ptm}{m}{n}
\rput(5.246875,-3.89){1}
\usefont{T1}{ptm}{m}{n}
\rput(6.446875,-2.69){1}
\usefont{T1}{ptm}{m}{n}
\rput(7.646875,-1.49){1}
\usefont{T1}{ptm}{m}{n}
\rput(8.846875,-0.29){1}
\usefont{T1}{ptm}{m}{n}
\rput(10.046875,0.91){1}
\usefont{T1}{ptm}{m}{n}
\rput(6.477031,-3.89){0}
\usefont{T1}{ptm}{m}{n}
\rput(7.6785936,-2.69){2}
\usefont{T1}{ptm}{m}{n}
\rput(8.877031,-1.49){0}
\usefont{T1}{ptm}{m}{n}
\rput(10.078594,-0.29){2}
\usefont{T1}{ptm}{m}{n}
\rput(7.646875,-3.89){1}
\usefont{T1}{ptm}{m}{n}
\rput(8.846875,-2.69){1}
\usefont{T1}{ptm}{m}{n}
\rput(10.046875,-1.49){1}
\usefont{T1}{ptm}{m}{n}
\rput(8.877031,-3.89){0}
\usefont{T1}{ptm}{m}{n}
\rput(10.077031,-2.69){0}
\usefont{T1}{ptm}{m}{n}
\rput(10.046875,-3.89){1}
\psline[linewidth=0.04cm,linestyle=dotted,dotsep=0.16cm](5.6,4.8)(5.6,5.6)
\psline[linewidth=0.04cm,linestyle=dotted,dotsep=0.16cm](5.6,-4.8)(5.6,-5.6)
\psline[linewidth=0.04cm,linestyle=dotted,dotsep=0.16cm](0.8,0.0)(0.0,0.0)
\psline[linewidth=0.04cm,linestyle=dotted,dotsep=0.16cm](10.6,0.0)(11.6,0.0)
\end{pspicture}
}
\caption{The activation process in $C_n\Box C_n$ when $4|n$.}
\end{figure}

\noindent If $m=n$ and $4|n$ then by the process shown in Figure 3 we may have equality in all above inequalities
and then $|D_0|=3n^2/8$. Then $wdyn(C_n\Box C_m)\leq |D_0|=3n^2/8$ when $4|n$.

\end{proof}


\section{A complexity result}

\noindent By the decision problem MINWDM we mean the problem of determining the smallest cardinality of any $\tau$-WDM in a given instance $(G,\tau)$. In this section we present a result concerning inapproximability of MINWDM. In proving our result we use a reduction from Minimum
Representative Problem (MINREP). In the following we explain the definition of MINREP from \cite{C}.

\noindent{\bf The Minimum Representative Problem (MINREP):}

\noindent Let $G = (A,B;E)$ be a given bipartite graph with the bipartite sets $A$ and $B$ and the edge set $E$. Assume that $A$ and $B$ are partitioned into
equal-sized subsets $A={\bigcup}_{i=1}^{\alpha}A_i$  and $B={\bigcup}_{j=1}^{\beta}B_j$. We have $|A_i|=|A|/\alpha$, for each $i$ and $|B_j|=|B|/\beta$, for each $j$. This partition of $V(G)$ induces a super-graph $H$ as follows.
Corresponding to each $A_i$ and $B_j$, $1\leq i\leq \alpha$, $1\leq j\leq \beta$, there are $\alpha + \beta$ super-vertices in $H$. There is a super-edge between $A_i$ and $B_j$ if there exist some
$a \in A_i$ and $b \in B_j$ such that $a$ and $b$ are adjacent in $G$.

\noindent We say a pair $(a , b)$ covers a super-edge $(A_i,B_j)$, if $a \in A_i$ and
$b \in B_j$ are adjacent in $G$.  We say $S \subseteq A_i\cup B_j$ covers a super-edge
$(A_i,B_j)$ if there exist $a , b \in S$ such that $(a , b)$ covers $(A_i,B_j)$.

\noindent The goal of the MINREP problem is to select the minimum number of representatives from each set $A_i$
and $B_j$ such that all super-edges are covered. That is, we wish to find subsets $A'\subseteq A$ and
$B'\subseteq B$ with the minimum total size $|A'|+|B'|$ such that, for every super-
edge $(A_i,B_j)$, there exist representatives $a\in A'\cap A_i$ and
$b\in B'\cap B_j$ that are adjacent in $G$. The following inapproximability result about MINREP was proved in \cite{R}. Note that $DTIME(n^{polylog(n)})$
is the class of problems which can be solved by a deterministic algorithm whose time complexity is bounded by
$n^{f(\log n)}$, where $f(\log n)$ is a polynomial in $\log n$.

\begin{thm}(\cite{R})
For any fixed $\epsilon > 0$, MINREP can not be approximated within the
ratio of ${\mathcal{O}}(2^{log^{(1-\epsilon)} n})$, unless $NP\subseteq DTIME(n^{polylog(n)})$.
\end{thm}


\noindent Inspired by the reduction technique of Chen \cite{C} (applied for dynamic monop-\
oly problem) we show the following result for MINWDM, whose proof is presented in Appendix A.

\begin{thm}\label{hardness}
For any fixed constant $\epsilon > 0$, MINWDM can not be approximated within
the ratio of ${\mathcal{O}}(2^{log^{1-\epsilon}n})$, unless
$NP\subseteq DTIME(n^{polylog(n)})$.
\end{thm}

\noindent In the following we give a remark concerning complexity of MINWDM in trees. In \cite{C}, a polynomial time algorithm for finding dynamic monopolies in trees with any threshold assignment is presented. This algorithm is based on the fact that any tree admits a dynamic monopoly with the smallest size which consists only of non-leaf vertices. But this property does not hold for weak dynamic monopolies of trees. Consider for example the tree $T$ with strict majority threshold depicted in the following figure, where each vertex (except some vertices of degree one) is identified by its label. Let
$D_0$ be any minimum WDM for $T$. It is easily seen that $D_0$ should contain $\{8,9,11,12,13,14,15\}$. Hence, for the corresponding sets $D_1$ and $D_2$ we have
$\{4,6,7\}\subseteq D_1$ and $3\in D_2$. Now if $10\notin D_0$, then either
$D_0=\{1,5,8,9,11,12,13,14,15\}$ or $D_0=\{1,2,8,9,11,12,13,14,15\}$. It follows that if $10\notin D_0$ then $|D_0|=9$.
But from other side $T$ contains a smaller WDM $\{8,9,11,12,13,14,15,10\}$ containing a vertex of degree one, i.e. 10.

\vspace{1cm}
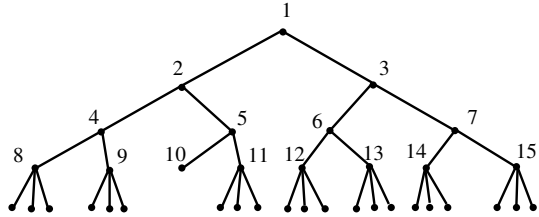
\begin{figure}[h]\label{tree-pic}
\centering
\scalebox{0.8} 
{
\begin{pspicture}(0,-1.7717187)(8.83625,1.7917187)
\psdots[dotsize=0.12](4.56,1.2482812)
\psdots[dotsize=0.12](2.88,0.32828125)
\psdots[dotsize=0.12](6.06,0.36828125)
\psdots[dotsize=0.12](1.54,-0.41171876)
\psdots[dotsize=0.12](3.72,-0.41171876)
\psdots[dotsize=0.12](5.34,-0.39171875)
\psdots[dotsize=0.12](7.42,-0.39171875)
\psdots[dotsize=0.12](0.44,-1.0117188)
\psdots[dotsize=0.12](1.68,-1.0517187)
\psdots[dotsize=0.12](2.88,-1.0117188)
\psdots[dotsize=0.12](3.86,-1.0117188)
\psdots[dotsize=0.12](4.88,-1.0117188)
\psdots[dotsize=0.12](6.0,-0.99171877)
\psdots[dotsize=0.12](6.94,-1.0117188)
\psdots[dotsize=0.12](8.44,-0.99171877)
\psline[linewidth=0.04cm](4.54,1.2682812)(2.9,0.36828125)
\psline[linewidth=0.04cm](2.9,0.36828125)(1.56,-0.39171875)
\psline[linewidth=0.04cm](1.56,-0.39171875)(0.46,-1.0117188)
\psline[linewidth=0.04cm](1.56,-0.43171874)(1.66,-1.0517187)
\psline[linewidth=0.04cm](3.72,-0.41171876)(2.9,-1.0117188)
\psline[linewidth=0.04cm](3.72,-0.37171876)(3.86,-1.0117188)
\psline[linewidth=0.04cm](3.7,-0.39171875)(2.86,0.36828125)
\psline[linewidth=0.04cm](4.54,1.2682812)(6.04,0.40828124)
\psline[linewidth=0.04cm](6.04,0.38828126)(5.34,-0.39171875)
\psline[linewidth=0.04cm](5.32,-0.41171876)(4.88,-0.99171877)
\psline[linewidth=0.04cm](5.32,-0.37171876)(6.0,-0.9717187)
\psline[linewidth=0.04cm](6.06,0.38828126)(7.42,-0.39171875)
\psline[linewidth=0.04cm](7.42,-0.39171875)(6.96,-0.99171877)
\psline[linewidth=0.04cm](7.42,-0.37171876)(8.46,-1.0117188)
\usefont{T1}{ptm}{m}{n}
\rput(4.626875,1.6182812){1}
\usefont{T1}{ptm}{m}{n}
\rput(2.8185937,0.65828127){2}
\usefont{T1}{ptm}{m}{n}
\rput(6.2476563,0.65828127){3}
\usefont{T1}{ptm}{m}{n}
\rput(1.4209375,-0.16171876){4}
\usefont{T1}{ptm}{m}{n}
\rput(3.8895311,-0.18171875){5}
\usefont{T1}{ptm}{m}{n}
\rput(5.1353126,-0.22171874){6}
\usefont{T1}{ptm}{m}{n}
\rput(7.7153125,-0.16171876){7}
\usefont{T1}{ptm}{m}{n}
\rput(0.16875,-0.8017188){8}
\usefont{T1}{ptm}{m}{n}
\rput(1.8951563,-0.82171875){9}
\usefont{T1}{ptm}{m}{n}
\rput(2.7725,-0.78171873){10}
\usefont{T1}{ptm}{m}{n}
\rput(4.156875,-0.8017188){11}
\usefont{T1}{ptm}{m}{n}
\rput(4.751406,-0.78171873){12}
\usefont{T1}{ptm}{m}{n}
\rput(6.063906,-0.76171875){13}
\usefont{T1}{ptm}{m}{n}
\rput(6.7726564,-0.7217187){14}
\usefont{T1}{ptm}{m}{n}
\rput(8.607031,-0.7217187){15}
\psdots[dotsize=0.12](0.06,-1.6717187)
\psdots[dotsize=0.12](0.38,-1.6917187)
\psdots[dotsize=0.12](0.68,-1.6917187)
\psdots[dotsize=0.12](1.38,-1.6717187)
\psdots[dotsize=0.12](1.68,-1.6917187)
\psdots[dotsize=0.12](1.92,-1.6917187)
\psline[linewidth=0.04cm](0.44,-1.0117188)(0.08,-1.6717187)
\psline[linewidth=0.04cm](0.44,-0.99171877)(0.4,-1.6917187)
\psline[linewidth=0.04cm](0.46,-1.0317187)(0.66,-1.6917187)
\psline[linewidth=0.04cm](1.66,-1.0517187)(1.38,-1.6717187)
\psline[linewidth=0.04cm](1.68,-1.0517187)(1.68,-1.6917187)
\psline[linewidth=0.04cm](1.68,-1.0117188)(1.92,-1.6917187)
\psdots[dotsize=0.12](3.52,-1.6717187)
\psdots[dotsize=0.12](3.8,-1.6717187)
\psdots[dotsize=0.12](4.14,-1.6717187)
\psdots[dotsize=0.12](4.64,-1.6917187)
\psdots[dotsize=0.12](4.92,-1.6917187)
\psdots[dotsize=0.12](5.26,-1.6917187)
\psdots[dotsize=0.12](5.78,-1.6917187)
\psdots[dotsize=0.12](6.08,-1.6717187)
\psdots[dotsize=0.12](6.36,-1.6717187)
\psdots[dotsize=0.12](6.76,-1.6717187)
\psdots[dotsize=0.12](7.0,-1.6717187)
\psdots[dotsize=0.12](7.3,-1.6717187)
\psdots[dotsize=0.12](8.12,-1.6717187)
\psdots[dotsize=0.12](8.44,-1.6717187)
\psdots[dotsize=0.12](8.72,-1.6717187)
\psline[linewidth=0.04cm](3.84,-0.99171877)(3.54,-1.6717187)
\psline[linewidth=0.04cm](3.86,-1.0317187)(3.8,-1.6717187)
\psline[linewidth=0.04cm](3.86,-0.99171877)(4.14,-1.6317188)
\psline[linewidth=0.04cm](4.88,-1.0117188)(4.64,-1.6717187)
\psline[linewidth=0.04cm](4.88,-0.99171877)(4.92,-1.6317188)
\psline[linewidth=0.04cm](4.9,-1.0317187)(5.24,-1.6917187)
\psline[linewidth=0.04cm](6.0,-0.9717187)(5.78,-1.6717187)
\psline[linewidth=0.04cm](6.0,-0.99171877)(6.08,-1.6117188)
\psline[linewidth=0.04cm](6.02,-0.9717187)(6.32,-1.6117188)
\psline[linewidth=0.04cm](6.92,-1.0317187)(6.78,-1.6717187)
\psline[linewidth=0.04cm](6.92,-0.99171877)(7.0,-1.5917188)
\psline[linewidth=0.04cm](6.94,-1.0317187)(7.28,-1.6117188)
\psline[linewidth=0.04cm](8.44,-0.9717187)(8.12,-1.6717187)
\psline[linewidth=0.04cm](8.46,-0.99171877)(8.44,-1.6317188)
\psline[linewidth=0.04cm](8.44,-0.95171875)(8.72,-1.6917187)
\end{pspicture}
}
\caption{A tree with a weak dynamic monopoly containing a leaf}
\end{figure}

\noindent We leave the determining of the smallest size of weak dynamic monopolies in trees as an unsolved question.

\section{Acknowledgment}

The authors thank the anonymous referee of the paper for his/her useful comments.

\section{Appendix A}

\noindent In the following we give the proof of Theorem \ref{hardness}.

\noindent {\bf Theorem.}
\noindent For any fixed constant $\epsilon > 0$, MINWDM can not be approximated within
the ratio of ${\mathcal{O}}(2^{log^{1-\epsilon}n})$, unless
$NP\subseteq DTIME(n^{polylog(n)})$.

\noindent\begin{proof} To prove this theorem we make a reduction from MINREP to
MINWDM. Let $G = (A,B;E)$ be an instance of MINREP with $N=|A|+|B|$ as the input size and let $M$  be the number of super-edges.
We use  basic gadgets $\Gamma_k$  shown in the following figure for the reduction. In fact,
$\Gamma_k$ is  $K_{2,k}$ in which each vertex $v_i$ ($i=1,..,k$) has $deg(v_i)=2$ and $\tau(v_i)=1$.

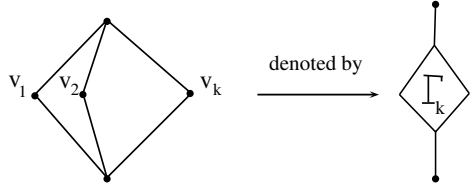
\begin{figure}[h]
\centering
\scalebox{0.8} 
{
\begin{pspicture}(0,-1.53)(7.66925,1.53)
\psdots[dotsize=0.12](1.63625,1.17)
\psdots[dotsize=0.12](1.63625,-1.45)
\psline[linewidth=0.027999999cm](1.61625,1.17)(0.43625,-0.05)
\psline[linewidth=0.027999999cm](0.43625,-0.05)(1.63625,-1.43)
\psline[linewidth=0.027999999cm](1.63625,-1.43)(1.23625,-0.03)
\psline[linewidth=0.027999999cm](1.23625,-0.03)(1.63625,1.19)
\psline[linewidth=0.027999999cm](1.63625,1.19)(3.03625,-0.03)
\psline[linewidth=0.027999999cm](3.03625,-0.03)(1.63625,-1.43)
\psline[linewidth=0.025999999cm](6.49625,-0.05)(7.07625,0.77)
\psline[linewidth=0.025999999cm](7.07625,0.77)(7.65625,-0.03)
\psline[linewidth=0.025999999cm](7.65625,-0.03)(7.09625,-0.67)
\psline[linewidth=0.025999999cm](7.09625,-0.67)(6.49625,-0.05)
\psline[linewidth=0.025999999cm](7.07625,0.77)(7.09625,1.45)
\psline[linewidth=0.025999999cm](7.09625,-0.65)(7.09625,-1.47)
\psline[linewidth=0.03cm,arrowsize=0.05291667cm 2.0,arrowlength=1.4,arrowinset=0.4]{->}(4.13625,-0.05)(6.15625,-0.05)
\usefont{T1}{ptm}{m}{n}
\rput(5.1092186,0.44){denoted by}
\psdots[dotsize=0.12](1.23625,-0.05)
\psdots[dotsize=0.12](3.01625,-0.03)
\psdots[dotsize=0.12](0.43625,-0.07)
\usefont{T1}{ptm}{m}{n}
\rput(0.10234375,0.165){\large v}
\usefont{T1}{ptm}{m}{n}
\rput(0.23796874,0.005){\footnotesize 1}
\usefont{T1}{ptm}{m}{n}
\rput(0.9423438,0.145){\large v}
\usefont{T1}{ptm}{m}{n}
\rput(1.10625,0.041){\footnotesize 2}
\usefont{T1}{ptm}{m}{n}
\rput(3.2823439,0.145){\large v}
\usefont{T1}{ptm}{m}{n}
\rput(3.4726562,0.045){\footnotesize k}
\psdots[dotsize=0.12](7.09625,1.45)
\psdots[dotsize=0.12](7.09625,-1.45)
\psline[linewidth=0.027999999cm](7.21625,0.25)(6.89625,0.27)
\psline[linewidth=0.027999999cm,tbarsize=0.07055555cm 5.0]{-|*}(6.97625,0.27)(6.97625,-0.19)
\usefont{T1}{ptm}{m}{n}
\rput(7.175,-0.27){\small k}
\psline[linewidth=0.027999999cm](7.19625,0.25)(7.19625,0.17)
\end{pspicture}
}
\caption{The basic gadget $\Gamma_k$.}
\end{figure}\label{gadget}

\noindent We construct now the graph $G'$ (corresponding to the instance $G$) as an instance of MINWDM.
The vertex set of $G'$ is partitioned into five subsets $V_1, V_2, \ldots, V_5$. The connection between these subsets is only via the gadgets
$\Gamma_k$ for appropriate values of $k$. In the following construction, for any two vertices $x,y\in V_1\cup \cdots \cup V_5$, when we say $x$ and $y$ are connected by a gadget $\Gamma_k$, it is simply meant that the two vertices $x$ and $y$ together with $k$ extra vertices form a gadget $\Gamma_k$, where $x$ and $y$ have degree $k$ in $\Gamma_k$. The details of the construction of $G'$ are as follow,

\begin{itemize}
\item $V_1=\{a|a \in A\}\cup\{b|b\in B\}$ and each vertex of $V_1$ has threshold $N^2$.

\item $V_2=\{u_{a,b}|(a,b)\in E\}$ and each vertex has threshold $2N^5$. Any vertex $u_{a,b}$ is connected to each of $a,b\in V_1$ by a basic gadget $\Gamma_{N^5}$.

\item $V_3=\{v_{i,j} | \text{$A_i$ is connected to $B_j$ by a super-edge} \}$ and each vertex has
      threshold $N^8$. If $a\in A_i$  and $b\in B_j$ then vertex $u_{a,b}\in V_2$ is connected
      to $v_{i,j}\in V_3$ by a basic gadget $\Gamma_{N^8}$.

\item $V_4=\{w_1,...,w_N\}$ and each vertex has threshold $M.N^2$. Each vertex $v_{i,j}\in V_3$
       is connected to each $w_k\in V_4$ by a basic gadget $\Gamma_{N^2}$, and each vertex
       $a,b\in V_1$ is  connected to each $w_k\in V_4$ by a basic gadget $\Gamma_{N}$.

\item  $V_5=\{z_1,...,z_N\}$ and each vertex has threshold $2.M.N^6$. Each vertex in $V_2$ is
       connected to each $z_k\in V_5$ by a gadget $\Gamma_{2N^4}$, and each vertex in $V_3$ is
       connected to each $z_k$ by a gadget $\Gamma_{2N^6}$.
\end{itemize}

\noindent The graph $G'$ is shown in Figure \ref{reduction}

\begin{figure}[h]
\centering
\includegraphics[scale=.6]{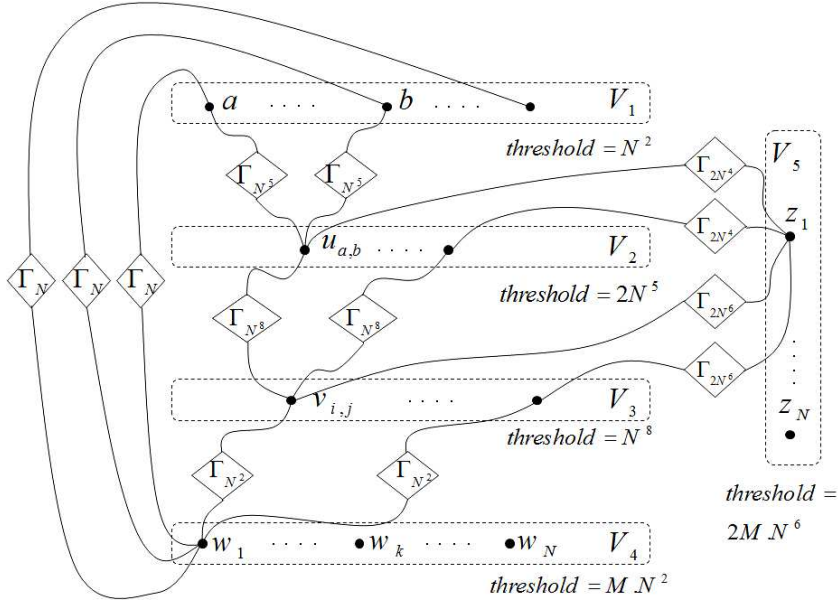}
\caption{The Graph $G'$.}\label{reduction}
\end{figure}

\noindent We show that the size of the optimal MINREP solution of $G$ is at most
twice the size of the  MINWDM of $G'$. Hence, any approximation algorithm for
MINWDM is transformed to an approximation algorithm for MINREP with the same performance ratio up to a constant factor.

\noindent Consider $A'\subseteq A$ and $B'\subseteq B$ as an optimal MINREP solution
of $G$. We show that $A'\cup B'\subseteq V_1$ is a WDM of $G'$.
Since $A'\cup B'$ is a MINREP solution, for any super-edge
$(A_i,B_j )$, there exist $a \in A' \cap A_i$ and $b \in B' \cap B_j$ such that
$\{a, b\} \in E$. The vertices of each gadget $\Gamma_{N^5}$ between $a$ or $b$ and $u_{a,b}$ are activated at time 1 and
hence, the vertex $u_{a,b}\in V_2$ is activated at time 2. Now the vertices of basic gadget $\Gamma_{N^8}$
between $u_{a,b}$ and $v_{i,j}$ become active at time 3 and then
$v_{i,j}\in V_3$ becomes active at time 4. This is true for
all super-edges, and thus all vertices in $V_3$ are active at time 4,
which implies that all vertices in $V_4$ and $V_5$ are active at time 6.
Therefore, all vertices in $V_1\setminus (A'\cup B')$ become active at time 8
by the vertices of $V_4$ and remaining vertices of $V_2$ become active at
time 8 by the vertices of $V_5$. Finally, in this manner all vertices
of $G'$ are activated.

\noindent Conversely, let $S$ be an optimal WDM of $G'$. It is obvious that
$|S|\leq N$ because $V_1$ is a WDM (as shown above). Note that $|V_2|<N^2$ and since $V_5$ is only connected to $V_2$ and $V_3$ (via gadgets), then the vertices of $V_5$ are only activated after the activation of all vertices of $V_3$. In fact the activation of whole $V_2$ and all but one vertex of $V_3$ does not imply the activation of $V_5$, since the threshold of each vertex in $V_5$ is $2MN^6$ and $2|V_2|N^4+2(M-
1)N^6 < 2MN^6$. Clearly the activation of the vertices of $V_3$ should be at the same time. Therefore either $V_3\subseteq S$ or
$V_3\cap S=\emptyset$. A similar argument shows that the vertices of $V_3$ become active only by the vertices of $V_2$.

\noindent In the following we show that we can modify $S$ such that $S\cap(V_4\cup V_5)=\emptyset$. The threshold of each vertex in $V_3$ is $N^8$. Hence none of the vertices in $V_3$ can be activated by the vertices of $V_4$. Now consider a vertex $a\in V_1$. Since $\tau(a)=N^2$, then $a$ is activated either by some vertex of $V_2$ which is connected to $a$ by $\Gamma_{N^5}$, or by the all vertices of $V_4$ which are connected to $a$ by $\Gamma_N$. In the second case $V_4\subseteq S$ and then $V_4=S$. Now we can replace $V_4$ by $V_1$ in $S$.

\noindent Similarly, no vertex of $V_3$ becomes active by the vertices of $V_5$. Consider $u_{a,b}\in V_2$ which is activated by some vertices of $V_5$. Note that in this case either $a$ or $b$ are not in $S$ because $\tau(u_{a,b})=2N^5$. Let $a\in S$ and $b\notin S$. As we mentioned before, the vertices of $V_5$ are only activated after the activation of all vertices of $V_3$, hence we can replace the vertices of $V_5$ in $S$ by $b$. Consequently, we may assume hereafter that $S\cap(V_4\cup V_5)=\emptyset$.

\noindent Also it is easily seen that the specified vertices $v_1,...,v_k$ from each gadget $\Gamma_k$ do not belong to $S$. The reason is simply that the minimum threshold in $V_1\cup\ldots\cup V_5$ is $N^2$ and $|S|\leq N$ by the optimality of $S$.

\noindent We may assume until far that $S\subset V_1\cup V_2\cup V_3$. We have also
$V_3\subseteq S$ or $V_3\cap S=\emptyset$. In the first case, we replace each $v_{i,j}\in V_3$ by $u_{a,b}\in V_2$, where $a\in A_i$ and $b\in B_j$. The resulting set still denoted by $S$ is WDM with the same cardinality of the previous $S$.
Finally, we replace any possible vertex $u_{a,b} \in S \cap V_2$ by the two vertices
$a$ and $b$ from $V_1$. Denote the resulting set by $S'$. Clearly $S'$ is a WDM and $|S'|\leq 2|S|$. Note that $S'\subseteq V_1$ and $S'$ is a (not necessarily optimal) solution for MINREP in the graph $G$. This completes the proof.
\end{proof}


\begin{thebibliography}{1}

\bibitem{ABW}
E. Ackerman, O. Ben-Zwi, G. Wolfovitz, Combinatorial Model and Bounds for Target Set Selection, Theoret. Comput. Sci. 411 (2010) 4017--4022.

\bibitem{ABST}
S.S. Adams, Z. Brass, C. Stokes, D.S. Troxell, Irreversible $k$-threshold and majority conversion processes on complete multipartite graphs and graph products, arXiv:1102.5361v1, 2011.

\bibitem{AM}
R.M. Anderson,  R.M. May, Infectious Diseases of Humans:
Dynamics and Control, Oxford University Press, 1991.

\bibitem{BM}
J.A. Bondy, U.S.R. Murty, Graph Theory, Springer 2008.

\bibitem{CL}
C-L. Chang, Y-D. Lyuu, On irreversible dynamic monopolies in general graphs, arXiv:0904.2306v3, 2009.

\bibitem{C}
N. Chen, On the approximability of influence in social networks,
SIAM J. Discrete Math. 23 (2009) 1400--1415.

\bibitem{CDPRS}
C.C. Centeno, M.C. Dourado, L.D. Penso, D. Rautenbach , J.L. Szwarcfiter, Irreversible conversion of graphs, Theoret. Comput. Sci. 412 (2011) 3693--3700.

\bibitem{DRi}
P. Domingos and M. Richardson, Mining the network value of customers, in: Proceedings of the 7th ACM International Conference on Knowledge Discovery and Data Mining, KDD, 2001, pp. 57--66.

\bibitem{DR}
P.A. Dreyer, F.S. Roberts, Irreversible $k$-threshold processes: Graph-theoretical threshold models of the spread of disease and of opinion, Discrete Appl. Math. 157 (2009) 1615--1627.

\bibitem{FFGHHKLS}
O. Favaron, G. Fricke, W. Goddard, S.M. Hedetniemi, S.T. Hedetniemi, P. Kristiansen, R.C. Laskar, R.D. Skaggs, Offensive alliances in graphs, Discuss. Math. Graph Theory 24 (2) (2004) 263--275.

\bibitem{F}
P. Flocchini, Contamination and decontamination in majority-based systems, Journal of Cellular Automata, 4(3):183-200, 2009.

\bibitem{FKRRS}
P. Flocchini, R. Kralovic, A. Roncato, P. Ruzicka, N. Santoro, On time versus size for monotone dynamic monopolies in regular topologies, J. Discrete Algorithms, 1 (2003) 129--150.

\bibitem{FLLPS}
P. Flocchini, E. Lodi, F. Luccio, L. Pagli, N. Santoro, Dynamic monopolies in tori, Discrete Appl. Math. 137 (2004) 197-–212.

\bibitem{KNSZ} K. Khoshkhah, M. Nemati, H. Soltani, M. Zaker, A study of monopolies in
graphs, Graphs and Combin. 29 (2013) 1417--1427.

\bibitem{KSZ} K. Khoshkhah, H. Soltani, M. Zaker, On dynamic monopolies of graphs: the average
and strict majority thresholds, Disc. Optimization 9 (2012) 77--83.

\bibitem{NNUW}
A. Nichterlein, R. Niedermeier, J. Uhlmann, M. Weller, On tractable cases of target set selection, Algorithms and computation. Part I, 378–-389,
Lecture Notes in Comput. Sci., 6506, Springer, Berlin, 2010.

\bibitem{R}
R. Raz, A Parallel Repetition Theorem, SIAM J. Computing, V.27(3), 763-803, 1998.

\bibitem{Z}
M. Zaker, On dynamic monopolies of graphs with general thresholds, Discrete Math. 312 (2012) 1136--1143.

\bibitem{Z2}
M. Zaker, Generalized degeneracy, dynamic monopolies and maximum degenerate subgraphs, Discrete Appl. Math., Discrete Appl. Math., 161 (2013) 2716–-2723.
\end{thebibliography}
\end{document}